\documentclass{amsart}
\usepackage[T1]{fontenc}
\usepackage[utf8]{inputenc}%pour les accents, sinon latin9
\usepackage{amsmath,amsthm,amsfonts,amscd,amssymb,eucal,latexsym,mathrsfs}
\usepackage{stmaryrd}
\usepackage{enumerate}
\usepackage{hyperref}
\usepackage[all]{xy}
\usepackage{etoolbox}

\usepackage{tikz,tikz-cd}
\usetikzlibrary{shapes.geometric}
\usetikzlibrary{arrows}

\usetikzlibrary{positioning}
\usepackage{float}
\usepackage{MnSymbol}
\usetikzlibrary{matrix}

\newcommand{\toc}{\tableofcontents}

\theoremstyle{plain}

\newtheorem*{theorem*}{Theorem}

\newtheorem*{corollary*}{Corollary}

\theoremstyle{definition}

\newtheorem*{definition*}{Definition}

\usetikzlibrary{calc,graphs}
\usepackage{xcolor}
\usetikzlibrary{arrows,decorations.pathmorphing}

\newcommand{\bb}[1]{\llbracket {#1}\rrbracket}
\newcommand{\stirl}[2]{\begin{Bmatrix}#1\\ #2\end{Bmatrix}}%stirling number second kind
\newcommand{\acts}{\curvearrowright}
\newcommand{\racts}{\curvearrowleft}
\newcommand{\G}{\Gamma}
\newcommand{\p}{\varphi}
\newcommand{\orb}{{\mathrm{orb}}}
\newcommand{\e}{\varepsilon}
\newcommand{\h}{\mathfrak{h}}
\newcommand{\pfk}{\mathfrak{p}}
\newcommand{\MG}{\mathbf{MG}}
\newcommand{\qfk}{\mathfrak{q}}
\newcommand{\f}{\mathrm{f}}
\newcommand{\s}{\mathbf{s}}
\renewcommand{\r}{\mathbf{r}}
\renewcommand{\d}{\mathrm{d}}
\newcommand{\opp}{\mathrm{opp}}
\newcommand{\IR}{\mathbb{R}}
\newcommand{\IH}{\mathbb{H}}
\newcommand{\CI}{\mathbb{C}}
\newcommand{\IF}{\mathbb{F}}
\newcommand{\IZ}{\mathbb{Z}}
\newcommand{\ZI}{\mathbb{Z}}
\newcommand{\IN}{\mathbb{N}}
\newcommand{\IQ}{\mathbb{Q}}
\newcommand{\QI}{\mathbb{Q}}
\newcommand{\SI}{\mathbb{S}}
\renewcommand{\Re}{\operatorname{Re}}

\DeclareMathOperator{\tr}{\mathrm{tr}}
\DeclareMathOperator{\ran}{\mathrm{ran}}
\DeclareMathOperator{\Prob}{\mathrm{Prob}}
\DeclareMathOperator{\dom}{\mathrm{dom}}
\DeclareMathOperator{\del}{\partial}
\DeclareMathOperator{\ssi}{\Leftrightarrow}
\DeclareMathOperator{\impl}{\Rightarrow}
\DeclareMathOperator{\inj}{\hookrightarrow}
\DeclareMathOperator{\linj}{\hookleftarrow}
\DeclareMathOperator{\limpl}{\Leftarrow}
\DeclareMathOperator{\lra}{\leftrightarrow}
\DeclareMathOperator{\la}{\leftarrow}
\DeclareMathOperator{\lar}{\leftarrow}
\DeclareMathOperator{\surj}{\twoheadrightarrow}
\DeclareMathOperator{\leftsurj}{\twoheadleftarrow}
\DeclareMathOperator{\lsurj}{\twoheadleftarrow}
\DeclareMathOperator{\id}{\mathrm{Id}}
\DeclareMathOperator{\hhom}{\mathrm{Hom}}
\DeclareMathOperator{\Con}{\mathrm{Con}}

\newcommand{\IE}{\mathbb{E}}

\newcommand{\m}{\mathfrak{m}}

\newcommand{\Ev}{\mathrm{Ev}}
\newcommand{\ev}{\mathrm{Ev}}
\newcommand{\ab}{\mathrm{ab}}

\newcommand{\SO}{\mathrm{SO}}
\newcommand{\SU}{\mathrm{SU}}

\DeclareMathOperator{\Frac}{\mathrm{Frac}}
\DeclareMathOperator{\Stab}{\mathrm{Stab}}
\DeclareMathOperator{\stab}{\mathrm{Stab}}
\DeclareMathOperator{\Tor}{Tor}
\DeclareMathOperator{\Inn}{Inn}
\DeclareMathOperator{\SAut}{\mathrm{SAut}}

\newcommand{\E}{\mathbf{E}}
\newcommand{\sgn}{\mathrm{sgn}}
\newcommand{\U}{U}
\newcommand{\TI}{\mathbb{T}}
\newcommand{\SL}{\mathrm{SL}}

\newcommand{\smatr}[4]{\left(\begin{smallmatrix} #1 & #2 \\ #3 & #4\end{smallmatrix}\right)}
\newcommand{\smatrix}[4]{\left(\begin{smallmatrix} #1 & #2 \\ #3 & #4\end{smallmatrix}\right)}
\newcommand{\matr}[4]{\left(\begin{array}{cc} #1 & #2 \\ #3 & #4\end{array}\right)}

\newcommand{\vect}[3]{\left(\begin{array}{c} #1 \\ #2 \\ #3 \end{array}\right)}

\newcommand{\conv}{\operatorname{conv}}

\newcommand{\ts}{\textsection}

\newcommand{\Act}{\mathrm{Act}}

\newcommand{\R}{{\mathbf R}}
\renewcommand{\r}{{\mathbf r}}
\renewcommand{\l}{{\mathbf l}}

\newcommand{\lb}{{\llbracket}}
\newcommand{\rb}{{\rrbracket}}

\newcommand{\ip}[1]{\langle#1\rangle} % \ip{a,b} gives us <a,b>
\newcommand{\norm}[1]{\|#1\|} 
\newcommand{\abs}[1]{|#1|} 

\newcommand{\dirlim}{\varinjlim}
\newcommand{\St}{\mathop{\mathrm{St}}}

\newcommand{\Aut}{\mathrm{Aut}}
\newcommand{\Id}{\mathrm{Id}}

\title{On the classification of CW complexes with prescribed links}
\author{Sylvain Barr\'e}
\author{Mika\"el Pichot}
\address{Sylvain Barr\'e, UMR 6205, LMBA, Université de Bretagne-Sud,BP 573, 56017, Vannes, France}\email{Sylvain.Barre@univ-ubs.fr}
\address{Mika\"el Pichot, McGill University, 805 Sherbrooke St W., Montr\'eal, QC H3A 0B9, Canada}\email{pichot@math.mcgill.ca}
\begin{document}

\setcounter{tocdepth}{1}

\begin{abstract}
Reporting on a computer--assisted search for nonpositively curved CW complexes of intermediate rank conducted some years ago. Not intended for publication.
\end{abstract}

\maketitle

 This is a brief report on a computer program, written between 2007 and 2009, aiming to classify certain CW complexes of dimension 2 with prescribed links, up to isomorphism of CW complex. Our initial motivation was the search of new groups of intermediate rank, which sometimes have specified link types, and the so-called ``flat closing conjecture''.  The program provided partial results, and was written to replace the more primitive version (based on a case-by-case analysis) presented in \cite[\ts 4]{rd}. 
 
 While our initial interest was on CW complexes having multiple vertices, the code was also tested in the single vertex case.  For example, our program provided the list of all  Moebius--Kantor  complexes having a single vertex. There are 27 such complexes in total. Here we call Möbius--Kantor complex, a CW complex of dimension 2 with triangle faces, whose links are isomorphic to the Moebius--Kantor graph.  The case-by-case analysis presented ealier in \cite[\ts 4]{rd} found a list 13 Moebius--Kantor complex on a single vertex,  with the additional property that all faces are positively oriented. 
 
The program tested for isomorphism of CW complexes. Thus, we claim, more precisely, that there exist exactly 27 Moebius--Kantor complexes having a single vertex, considered up to isomorphism of CW-complexes.  We do not  have a formal proof of this fact: the list of complexes was returned by the computer, and a ``human readable'' proof of completeness (which we have not attempted to write) would quite likely occupy pages of unpleasant computations.  The ``complete'' list of these complexes is shown in \ts\ref{S - 27 MK complexes}. We note, in addition, that the use of a computer program to search for group presentations is quite standard. In \cite{CMSZ},  Cartwright, Mantero, Steger, and Zappa  used a similar computer theoretic  approach to find all ``triangle presentations'' of groups acting on buildings in the $\tilde A_2$ case. 

Our program was mostly used to find complexes with multiple vertices. Unfortunately, it suffered from its (relatively speaking) ``success'', as it  provided long lists of CW complexes which remained mostly unstudied. For example, the program provided a list of more than 200 Moebius--Kantor complexes on three vertices: a list which is in fact incomplete, since the computation was aborted after a few days of computations.

   One of our initial interests in \cite{rd} was the ``flat closing conjecture''.  Let  $G$  be a countable acting geometrically on a CAT(0) complex $X$ of dimension 2. The conjecture states that if $X$ contains a flat plane $(\simeq \IR^2$), then $G$ contains $\ZI^2$. This conjecture was known in two situations:
   \begin{enumerate}
   \item  $G$ is (relatively) hyperbolic  
   \item $G$ acts on a Euclidean buildings. 
   \end{enumerate}
Both cases of (relatively) hyperbolic groups, and of groups acting on Euclidean buildings, are thoroughly studied in geometric group theory. Since ref.\ \cite{rd} we have referred, informally, to groups  sitting  ``in between'' these two cases, as groups of intermediate rank. The flat closing conjecture remains open for these groups in general. A formal property was also introduced in \cite{rd}, called the ``exponential mesoscopic rank'' property, to insure that $X$  is ``of intermediate rank'' in the sense that it is neither a Euclidean building, nor a space with isolated flats.

There are no counter-example to the closing conjecture among the list of 13 found in \cite{rd}, and we were not able to find a counterexample in larger lists found afterwards, despite some efforts. We have  since isolated  two particularly interesting simply connected Moebius--Kantor complexes of intermediate rank: the even and the odd Moebius--Kantor complex (both have the mesoscopic rank property). Of course, many classical groups were already known to be ``of intermediate rank''. A typical example is $\Aut(F_2)$. Yet Moebius--Kantor complexes seemed particularly attractive, and not that far, at least locally, from being Euclidean buildings.

The program was also tested on Euclidean buildings themselves, and we found a rather nice  example of a $\tilde B_2$ building on 7 vertices having 45 faces. In fact, the program  found, in this case too, hundreds of such complexes, but we did not have the energy to test for isomorphism. One example of such a complex is reproduced in \ts \ref{S - tilde b2}.  

We shall describe very quickly how the program was implemented in \ts\ref{S - code}. Ultimately, the computer approach to intermediate rank geometry was not as successful as we had hoped for, and it was essentially abandoned in 2009. We have since focused  on more abstract approaches to intermediate rank geometry not relying on the use of computers.
We thank Alex Lou\'e for private communications on his (independent) computer searches for Moebius--Kantor complexes, which are at the origin of the present report. 
 
\section{The 27 Moebius--Kantor complexes on a single vertex}\label{S - 27 MK complexes}
 
 The following is a list of Moebius Kantor complexes on a single vertex, as returned by our program. An explanation for this list is 
 provided in \ts\ref{S - code}.
 
 $V_{1}$: 
\[
\begin{vmatrix}
(1,0) &(10,14) &(15,8) &(13,11) &(15,8) &(3,7) &(6,4) &(2,5)\\
(5,2) &(3,7) &(9,12) &(12,9) &(13,11) &(6,4) &(9,12) &(6,4)\\
(1,0) &(8,15) &(1,0) &(10,14) &(10,14) &(5,2) &(13,11) &(3,7)\\
\end{vmatrix}
\] 

 $V_{2}$: 
\[
\begin{vmatrix}
(1,0) &(4,3) &(7,6) &(13,8) &(15,10) &(15,10) &(7,6) &(12,9)\\
(5,2) &(2,5) &(5,2) &(9,12) &(11,14) &(9,12) &(11,14) &(4,3)\\
(1,0) &(4,3) &(7,6) &(13,8) &(15,10) &(1,0) &(13,8) &(14,11)\\
\end{vmatrix}
\]

 $V_{3}$: \[
\begin{vmatrix}
(1,11) &(15,4) &(5,2) &(2,5) &(12,6) &(8,3) &(4,15) &(10,7)\\
(6,12) &(9,14) &(3,8) &(6,12) &(7,10) &(4,15) &(10,7) &(8,3)\\
(13,0) &(13,0) &(13,0) &(11,1) &(11,1) &(14,9) &(2,5) &(14,9)\\
\end{vmatrix}
\]

 $V_{4}$: \[
\begin{vmatrix}
(1,11) &(15,6) &(5,14) &(2,9) &(9,2) &(12,10) &(3,8) &(7,4)\\
(10,12) &(5,14) &(3,8) &(10,12) &(7,4) &(15,6) &(7,4) &(5,14)\\
(13,0) &(13,0) &(13,0) &(11,1) &(3,8) &(11,1) &(9,2) &(15,6)\\
\end{vmatrix}
\]

 $V_{5}$: \[
\begin{vmatrix}
(1,0) &(4,3) &(10,13) &(9,8) &(12,11) &(15,6) &(7,14) &(5,2)\\
(5,2) &(2,5) &(14,7) &(13,10) &(10,13) &(11,12) &(3,4) &(7,14)\\
(1,0) &(4,3) &(6,15) &(9,8) &(12,11) &(1,0) &(9,8) &(15,6)\\
\end{vmatrix}
\]

 $V_{6}$: \[
\begin{vmatrix}
(1,0) &(4,3) &(10,9) &(15,8) &(8,15) &(10,9) &(5,2) &(6,13)\\
(5,2) &(2,5) &(8,15) &(7,12) &(14,11) &(4,3) &(7,12) &(14,11)\\
(1,0) &(4,3) &(10,9) &(1,0) &(6,13) &(14,11) &(13,6) &(12,7)\\
\end{vmatrix}
\]

 $V_{7}$: \[
\begin{vmatrix}
(1,0) &(4,3) &(15,10) &(15,10) &(6,8) &(13,7) &(13,7) &(12,9)\\
(5,2) &(2,5) &(11,14) &(9,12) &(9,12) &(8,6) &(2,5) &(4,3)\\
(1,0) &(4,3) &(15,10) &(1,0) &(13,7) &(11,14) &(6,8) &(14,11)\\
\end{vmatrix}
\]

 $V_{8}$: \[
\begin{vmatrix}
(13,12) &(1,14) &(15,6) &(5,8) &(2,4) &(3,9) &(8,5) &(3,9)\\
(1,14) &(15,6) &(7,10) &(13,12) &(9,3) &(8,5) &(6,15) &(10,7)\\
(13,12) &(11,0) &(11,0) &(11,0) &(14,1) &(4,2) &(10,7) &(2,4)\\
\end{vmatrix}
\]

 $V_{9}$: \[
\begin{vmatrix}
(1,14) &(15,8) &(5,12) &(2,11) &(14,1) &(3,10) &(9,6) &(10,3)\\
(15,8) &(9,6) &(11,2) &(10,3) &(12,5) &(9,6) &(5,12) &(4,13)\\
(7,0) &(7,0) &(7,0) &(14,1) &(4,13) &(11,2) &(13,4) &(8,15)\\
\end{vmatrix}
\]

 $V_{10}$: \[
\begin{vmatrix}
(1,0) &(4,3) &(15,7) &(12,8) &(6,14) &(10,13) &(10,13) &(5,2)\\
(5,2) &(2,5) &(8,12) &(9,11) &(3,4) &(8,12) &(14,6) &(7,15)\\
(1,0) &(4,3) &(1,0) &(10,13) &(9,11) &(11,9) &(7,15) &(14,6)\\
\end{vmatrix}
\]

 $V_{11}$: \[
\begin{vmatrix}
(1,0) &(14,3) &(9,8) &(11,10) &(15,6) &(12,7) &(4,13) &(3,14)\\
(5,2) &(4,13) &(13,4) &(15,6) &(7,12) &(8,9) &(12,7) &(15,6)\\
(1,0) &(14,3) &(9,8) &(11,10) &(1,0) &(10,11) &(2,5) &(5,2)\\
\end{vmatrix}
\]

 $V_{12}$: \[
\begin{vmatrix}
(1,0) &(14,3) &(15,6) &(11,9) &(4,13) &(11,9) &(4,13) &(3,14)\\
(5,2) &(4,13) &(7,12) &(10,8) &(8,10) &(8,10) &(12,7) &(15,6)\\
(1,0) &(14,3) &(1,0) &(7,12) &(11,9) &(15,6) &(2,5) &(5,2)\\
\end{vmatrix}
\]

 $V_{13}$: \[
\begin{vmatrix}
(1,0) &(14,3) &(15,4) &(6,13) &(3,14) &(10,7) &(9,12) &(4,15)\\
(5,2) &(4,15) &(9,12) &(12,9) &(13,6) &(6,13) &(11,8) &(10,7)\\
(1,0) &(14,3) &(1,0) &(8,11) &(5,2) &(8,11) &(7,10) &(2,5)\\
\end{vmatrix}
\]

 $V_{14}$: \[
\begin{vmatrix}
(14,13) &(1,9) &(15,7) &(5,2) &(14,13) &(2,5) &(10,4) &(12,3)\\
(12,3) &(4,10) &(8,6) &(3,12) &(8,6) &(6,8) &(5,2) &(4,10)\\
(14,13) &(11,0) &(11,0) &(11,0) &(7,15) &(9,1) &(7,15) &(9,1)\\
\end{vmatrix}
\]

 $V_{15}$: \[
\begin{vmatrix}
(13,12) &(1,14) &(15,6) &(5,8) &(2,9) &(10,3) &(3,10) &(12,13)\\
(1,14) &(15,6) &(5,8) &(9,2) &(10,3) &(4,11) &(11,4) &(8,5)\\
(13,12) &(7,0) &(7,0) &(7,0) &(14,1) &(6,15) &(9,2) &(4,11)\\
\end{vmatrix}
\]

 $V_{16}$: \[
\begin{vmatrix}
(6,5) &(10,9) &(1,8) &(15,12) &(5,6) &(2,13) &(14,7) &(12,15)\\
(4,11) &(4,11) &(13,2) &(11,4) &(7,14) &(14,7) &(2,13) &(10,9)\\
(6,5) &(10,9) &(3,0) &(3,0) &(3,0) &(8,1) &(12,15) &(8,1)\\
\end{vmatrix}
\]

 $V_{17}$: \[
\begin{vmatrix}
(8,7) &(13,12) &(1,14) &(15,5) &(5,15) &(2,11) &(11,2) &(3,10)\\
(6,9) &(1,14) &(15,5) &(6,9) &(10,3) &(10,3) &(7,8) &(9,6)\\
(8,7) &(13,12) &(4,0) &(4,0) &(4,0) &(14,1) &(13,12) &(11,2)\\
\end{vmatrix}
\]

 $V_{18}$: \[
\begin{vmatrix}
(1,14) &(15,8) &(5,2) &(2,5) &(12,7) &(3,10) &(9,6) &(10,3)\\
(3,10) &(9,6) &(7,12) &(4,13) &(8,15) &(9,6) &(7,12) &(4,13)\\
(11,0) &(11,0) &(11,0) &(14,1) &(14,1) &(5,2) &(13,4) &(8,15)\\
\end{vmatrix}
\]

 $V_{19}$: \[
\begin{vmatrix}
(12,1) &(1,12) &(15,3) &(5,7) &(14,10) &(2,13) &(15,3) &(7,5)\\
(2,13) &(11,4) &(14,10) &(6,8) &(11,4) &(8,6) &(2,13) &(4,11)\\
(12,1) &(9,0) &(9,0) &(9,0) &(3,15) &(5,7) &(14,10) &(6,8)\\
\end{vmatrix}
\]

 $V_{20}$: \[
\begin{vmatrix}
(15,10) &(12,11) &(1,14) &(5,2) &(2,5) &(12,11) &(3,8) &(13,6)\\
(9,0) &(6,13) &(3,8) &(7,4) &(6,13) &(10,15) &(7,4) &(7,4)\\
(15,10) &(12,11) &(9,0) &(9,0) &(14,1) &(14,1) &(5,2) &(3,8)\\
\end{vmatrix}
\]

 $V_{21}$: \[
\begin{vmatrix}
(15,10) &(13,12) &(1,3) &(15,10) &(5,9) &(2,8) &(9,5) &(12,13)\\
(11,14) &(11,14) &(14,11) &(9,5) &(4,7) &(7,4) &(4,7) &(8,2)\\
(15,10) &(13,12) &(6,0) &(6,0) &(6,0) &(3,1) &(2,8) &(3,1)\\
\end{vmatrix}
\]

 $V_{22}$: \[
\begin{vmatrix}
(11,10) &(1,13) &(15,2) &(5,7) &(2,15) &(13,1) &(3,6) &(7,5)\\
(9,12) &(12,9) &(7,5) &(6,3) &(14,8) &(12,9) &(11,10) &(6,3)\\
(11,10) &(4,0) &(4,0) &(4,0) &(13,1) &(8,14) &(15,2) &(14,8)\\
\end{vmatrix}
\]

 $V_{23}$: \[
\begin{vmatrix}
(9,8) &(15,10) &(1,6) &(15,10) &(5,2) &(12,3) &(3,12) &(14,11)\\
(13,4) &(11,14) &(5,2) &(9,8) &(1,6) &(14,11) &(13,4) &(12,3)\\
(9,8) &(15,10) &(7,0) &(7,0) &(7,0) &(6,1) &(5,2) &(4,13)\\
\end{vmatrix}
\]

 $V_{24}$: \[
\begin{vmatrix}
(1,5) &(15,12) &(5,1) &(15,12) &(11,2) &(10,3) &(7,4) &(11,2)\\
(4,7) &(1,5) &(2,11) &(13,9) &(3,10) &(14,8) &(9,13) &(7,4)\\
(6,0) &(6,0) &(6,0) &(8,14) &(15,12) &(13,9) &(14,8) &(3,10)\\
\end{vmatrix}
\]

 $V_{25}$: \[
\begin{vmatrix}
(1,6) &(15,12) &(5,10) &(2,13) &(13,2) &(15,12) &(12,15) &(5,10)\\
(7,4) &(13,2) &(9,14) &(8,11) &(7,4) &(11,8) &(10,5) &(11,8)\\
(3,0) &(3,0) &(3,0) &(6,1) &(9,14) &(9,14) &(6,1) &(7,4)\\
\end{vmatrix}
\]

 $V_{26}$: \[
\begin{vmatrix}
(1,6) &(15,13) &(5,7) &(15,13) &(14,9) &(12,3) &(3,12) &(11,4)\\
(5,7) &(12,3) &(6,1) &(14,9) &(10,8) &(4,11) &(11,4) &(5,7)\\
(2,0) &(2,0) &(2,0) &(8,10) &(13,15) &(6,1) &(9,14) &(8,10)\\
\end{vmatrix}
\]

 $V_{27}$: \[
\begin{vmatrix}
(1,8) &(15,2) &(5,12) &(2,15) &(14,11) &(11,14) &(14,11) &(6,3)\\
(13,4) &(7,10) &(1,8) &(10,7) &(6,3) &(3,6) &(12,5) &(4,13)\\
(9,0) &(9,0) &(9,0) &(8,1) &(2,15) &(7,10) &(4,13) &(12,5)\\
\end{vmatrix}
\]

\section{Code description}\label{S - code}

The computer program was elementary and restricted to triangle faces, although it would not be difficult to include more general faces. Furthermore, the code works for general links not necessarily related to ``intermediate rank''.  

The first step is to record the desired links as graphs, with an injective numbering for the link vertices.  This can also be viewed as a single graph with marked sets of edges. In this case, a vertex of the would-be CW complex is represented as a set of edges corresponding to a given marking. 

The second step is to create the ``potential faces''. This is a list of all faces which are allowed in the CW complexes (any complex) we are building.  For instance, one can create here all possible faces that visit three given links, or all faces having two vertices in one link and the other vertex in a second link.  

Individually, each potential face is represented as a triple of pairs, where a given pair (called a ``wire'') represents a face edge from a vertex link to another. For example, the first face of $V_1$ in \ts\ref{S - 27 MK complexes} is $(1,0), (5,2),(1,0)$, which represent a triangle face spanning the vertices $0,1,2,5$ in the link, where the link is the Moebius--Kantor graph and its vertices are numbered cyclically from 0 to 15. We call such faces, i.e., triples of pairs, ``wired faces'', and  a set of links together with a set of pairwise compatible wired faces, a ``partial wired complex''. The pairwise compatibility test is straightforward: the wires should either be equal, inverse of each other, or disjoint, for two faces to be compatible in a potential CW complex.  

In practice, the list of potential faces is a list of triples of pairs having a \emph{very large} number of elements. This list will be trimmed incrementally as the wired complexes are  built recursively.

It is of course straightforward to translate a wired complex into a standard presentation complex. For example, say the wired complex (viewed as a list of wired faces) on the Moebius--Kantor link, with vertices numbered from 1 to 16, is:
\begin{align*}
&((1,2)(1,2)(3,16))\\
&((5,6)(1,2)(13,4))\\
&((3,16)(15,12)(13,4))\\
&((3,16)(11,14)(9,8))\\
&((5,6)(5,6)(7,10))\\
&((11,14)(13,4)(15,12))\\
&((11,14)(15,12)(7,10))\\
&((9,8)(9,8)(7,10))
\end{align*}
then the corresponding presentation complex is:
\[
[[1,1,2], [3,1,4], [2,5,4], [2,6,7], [3,3,8], [6,4,5], [6,5,8], [7,7,8]].
\]
which is one of the complexes listed in \cite[\ts 4]{rd}.

    The description of CW complexes through their links, as wired complexes, originates in the theory of graph of groups. The CW complexes thus constructed are polygonal $L$-complexes, where $L$ is the set of admissible links (see \cite{BB}). Our description with wires can easily be formalized using an order 2 bijection in place of wires (compare for instance \cite[\ts 3]{Ben}).

In the third step, the program debuts with a partial wired complex (for example, one can give a single wired face, or two disjoint faces, or even, for very small links, no initial face at all) and runs through the list of potential faces in order to determine which faces are compatible with the partially defined complex. It creates a list of newly compatible faces, adds them one by one to a potential CW complex, and recurses. At some point, the list of newly compatible faces will be empty: when it is empty, we have obtained a new complex if and only if every link edge (equivalently, link vertex) has been wired. If not, the complex found is  a ``partial wired complex'', and is discarded. To summarize, the third step is an elementary recursive function which, given a partially constructed complex, computes first the faces which are compatible with the complex, and then recurse on these faces augmenting it one face at a time.  

Finally, the program tests for isomorphism using the automorphism group of a graphs as input.  The use of wired complex is particularly adapted to this task. In our implementation, the test was time consuming, especially since large lists of complexes were typically discovered by the program (even for a single vertex complex) and every newly discovered complex was tested against all previously discovered examples (which can be in the hundreds). A more memory intensive approach would record every discovered examples and test for isomorphism in a final step. We have also tried to ``break the symmetry'' by entering well--chosen sets of initial faces, or by limiting the set of created faces to visit certain types of vertices in certain ways, in order to find smaller list (but less canonical) of complexes.

The program was written  in C++.

\section{An example of a complex of type $\tilde B_2$}\label{S - tilde b2}

Here is how the $\tilde B_2$ example was returned by our program:

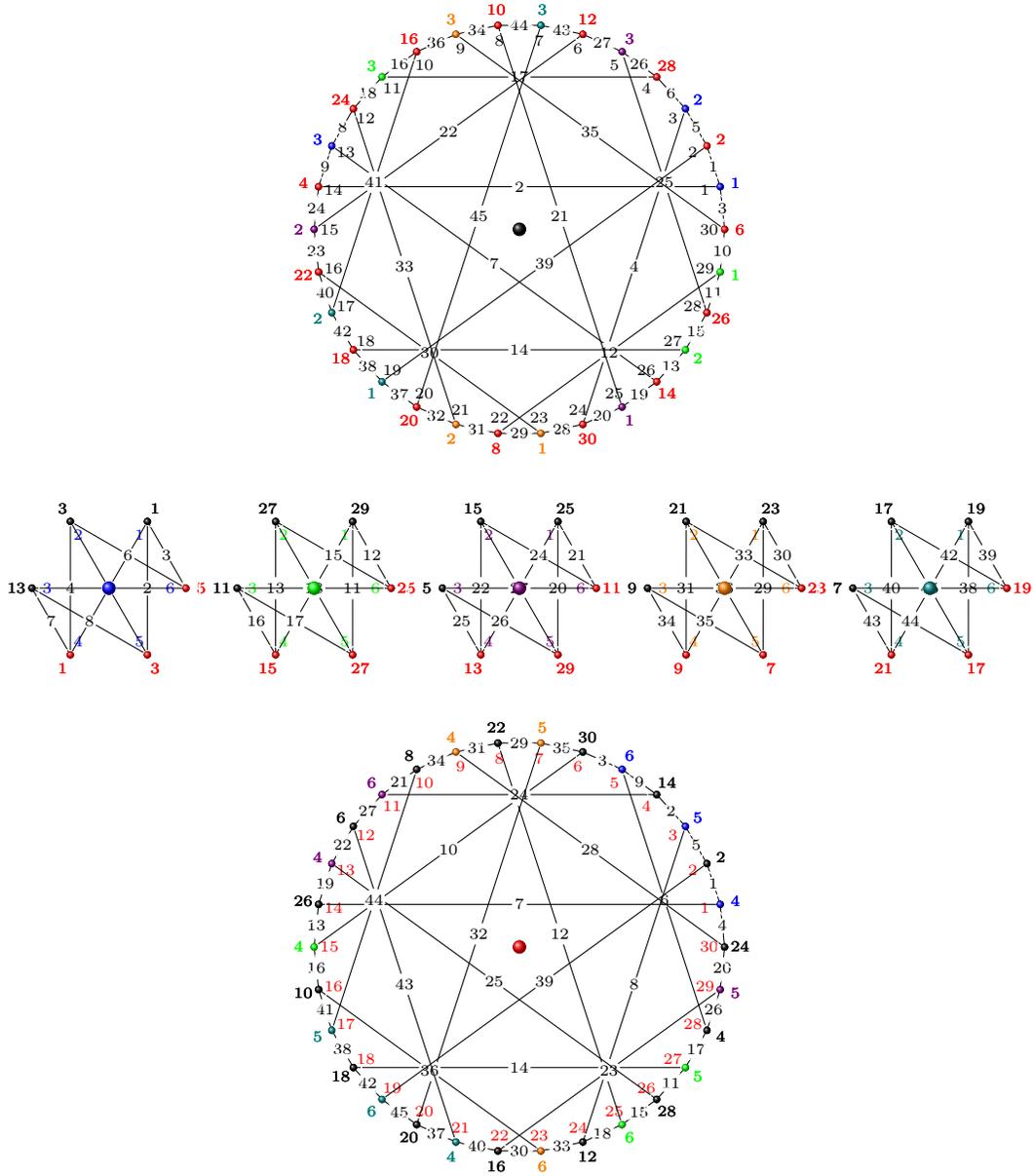
\begin{figure}[H]
\centering
\centerline{
}
\end{figure}

This is a complex with 7 vertices (corresponding to the 7 links above)
 having 45 faces. The links were chosen to insure the universal cover of this complex is Euclidean a building of type $\tilde B_2$ (a property that can be verified locally). The numbers inside the drawing above represent the wired faces, while the standard presentation CW complex (with 7 vertices drawn in various colours) is given by

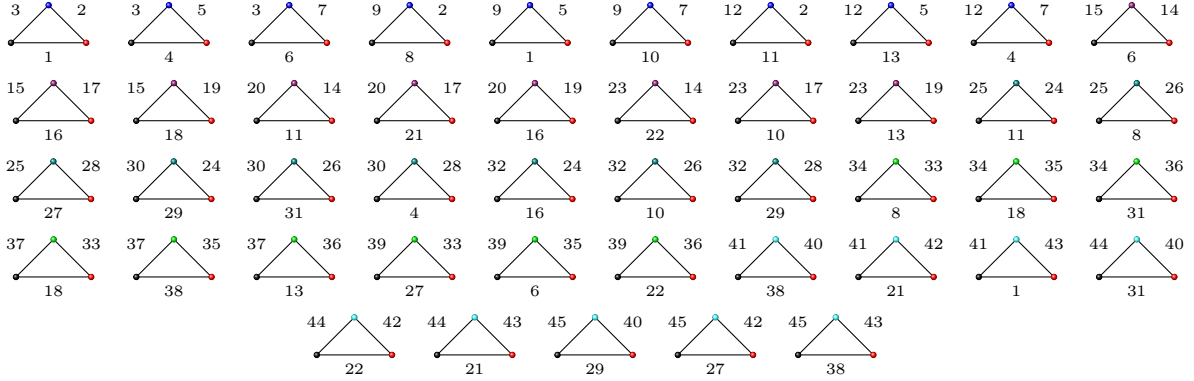
\begin{figure}[H]
\centering
\centerline{%
}
\end{figure}

In the case of the $\tilde B_2$ complex, the program was run on a more powerful computer  provided by the university of Tokyo IPMU (Institute for the Physics and Mathematics of the Universe) at Kashiwa (and aimed, in principle, for celestial computations! This is an opportunity for us to thank the university's IT staff for their constant technical support with the computers, and for providing us unrestricted access).  While the $\tilde B_2$ example has 7 vertices, which may seem quite large, an obvious choice for the potential wired faces was to list first the potential faces that fully complete only one of the five $K_{3,3}$ graph first. This gave us bricks, in which one of the five $K_{3,3}$ graph is complete. Five different types of bricks could then be assembled together in order to create the above complex. In fact, this complex came together with hundreds of  siblings, which we did not test for isomorphisms (the automorphism groups of the above graphs being quite large). We do not know how many such complexes exist up to isomorphism, although this seems achievable with additional work.

While these $\tilde B_2$ complexes remained unpublished, they lead to a new idea to construct groups of intermediate rank: that of removing chambers at random starting from a building. For example, one can remove one of the 45 faces in the $\tilde B_2$ example above at random, and obtain a new group acting on a new space which ``ought to have many flats''. We began to explore this idea in \cite{chambers}. The main difficult to find groups of intermediate rank is ensure that flats do not disappear in the process, and having a large ``local rank'' (which is the case when one remove one edge from a spherical building) is certainly insufficient to achieve this (unless, of course, the local rank is maximal, i.e., 2).  For example, there were already spaces with the isolated flat property present in the list of 13 examples from \cite{rd}.

\section{On the construction of large examples}

Upon  inspection of various folders, it seems that the largest example of a Moebius--Kantor complex we had considered had a total of 24 vertices and the following list of 192 faces: 

\bigskip

\noindent 
[[0,20,10], [1,720,10], [1,21,11], [2,721,11], [2,22,12], [3,722,12], [3,23,13], [4,723,13], [4,24,14], [5,724,14], [5,25,15], [6,725,15], [6,26,16], [7,726,16], [7,27,17], [0,727,17], [100,120,110], [101,20,110], [101,121,111], [102,21,111], [102,122,112], [103,22,112], [103,123,113], [104,23,113], [104,124,114], [105,24,114], [105,125,115], [106,25,115], [106,126,116], [107,26,116], [107,127,117], [100,27,117], [200,220,210], [201,120,210], [201,221,211], [202,121,211], [202,222,212], [203,122,212], [203,223,213], [204,123,213], [204,224,214], [205,124,214], [205,225,215], [206,125,215], [206,226,216], [207,126,216], [207,227,217], [200,127,217], [300,320,310], [301,220,310], [301,321,311], [302,221,311], [302,322,312], [303,222,312], [303,323,313], [304,223,313], [304,324,314], [305,224,314], [305,325,315], [306,225,315], [306,326,316], [307,226,316], [307,327,317], [300,227,317], [400,420,410], [401,320,410], [401,421,411], [402,321,411], [402,422,412], [403,322,412], [403,423,413], [404,323,413], [404,424,414], [405,324,414], [405,425,415], [406,325,415], [406,426,416], [407,326,416], [407,427,417], [400,327,417], [500,520,510], [501,420,510], [501,521,511], [502,421,511], [502,522,512], [503,422,512], [503,523,513], [504,423,513], [504,524,514], [505,424,514], [505,525,515], [506,425,515], [506,526,516], [507,426,516], [507,527,517], [500,427,517], [600,620,610], [601,520,610], [601,621,611], [602,521,611], [602,622,612], [603,522,612], [603,623,613], [604,523,613], [604,624,614], [605,524,614], [605,625,615], [606,525,615], [606,626,616], [607,526,616], [607,627,617], [600,527,617], [700,720,710], [701,620,710], [701,721,711], [702,621,711], [702,722,712], [703,622,712], [703,723,713], [704,623,713], [704,724,714], [705,624,714], [705,725,715], [706,625,715], [706,726,716], [707,626,716], [707,727,717], [700,627,717], [10,6,220], [11,7,221], [12,0,222], [13,1,223], [14,2,224], [15,3,225], [16,4,226], [17,5,227], [110,106,320], [111,107,321], [112,100,322], [113,101,323], [114,102,324], [115,103,325], [116,104,326], [117,105,327], [210,206,420], [211,207,421], [212,200,422], [213,201,423], [214,202,424], [215,203,425], [216,204,426], [217,205,427], [310,306,520], [311,307,521], [312,300,522], [313,301,523], [314,302,524], [315,303,525], [316,304,526], [317,305,527], [410,406,620], [411,407,621], [412,400,622], [413,401,623], [414,402,624], [415,403,625], [416,404,626], [417,405,627], [510,506,720], [511,507,721], [512,500,722], [513,501,723], [514,502,724], [515,503,725], [516,504,726], [517,505,727], [610,606,20], [611,607,21], [612,600,22], [613,601,23], [614,602,24], [615,603,25], [616,604,26], [617,605,27], [710,706,120], [711,707,121], [712,700,122], [713,701,123], [714,702,124], [715,703,125], [716,704,126], [717,705,127]]

\bigskip

This raises the question of constructing, especially by hand, infinitely many (arbitrary large) compact Moebius--Kantor complexes. This question was answered later in \cite{surgery} by using a different method. We do not know how to estimate the number of Moebius--Kantor complexes (or other types of complexes of intermediate rank) on $n$ vertices.


\begin{thebibliography}{00}


\bibitem{BB} Ballmann W., Brin M., Polygonal complexes and combinatorial group theory, Geom. Dedicata, 50 (!994), 165--191.


\bibitem{rd} Barr\'e, S.,  Pichot, M., Intermediate rank and property RD, preprint, 2007.

\bibitem{surgery} Barr\'e, S.,  Pichot, M., Surgery on discrete groups, preprint, 2016.

\bibitem{chambers}
Barr{\'e} and Pichot.
\newblock Removing chambers in Bruhat--Tits buildings.
\newblock {\em Israel J. Math.}, 202:117--160 (2014).

\bibitem{Ben} Benakli N, Polygonal complexes I: Combinatorial and geometric properties, J. Pure and Applied Alg., 97, (1994) 247--263.


\bibitem{CMSZ} Cartwright, Donald I.; Mantero, Anna Maria; Steger, Tim; Zappa, Anna Groups acting simply transitively on the vertices of a building of type $\tilde A\sb 2$. I \& II. The cases $q=2$ and $q=3$. Geom. Dedicata 47 (1993), no. 2, 143--166 \& 167--223.

\end{thebibliography}
\end{document}